\documentclass[12pt]{article}
\usepackage{amsmath,amsfonts}
\usepackage{times}
\usepackage{epstopdf,xcolor,cite,epstopdf}
\usepackage{caption}
\usepackage{hyperref}
\hypersetup{colorlinks=true,linkcolor=blue,citecolor=red}
\usepackage{tikz}
\usetikzlibrary{calc}
\usetikzlibrary{trees}
\tikzstyle{every node}=[circle,inner sep=1pt,fill=white!60]
\tikzstyle{tn}=[shape=circle, draw, color=black!70]
\tikzstyle{tn1}=[shape=circle, draw]
\tikzstyle{tn2}=[shape=circle, draw,inner sep=1.5pt]
\tikzstyle{tn3}=[shape=rectangle, draw,inner sep=1.5pt]
\tikzstyle{marke}=[shape=circle,minimum size=0.1cm, draw,blue]

\newtheorem{thm}{Theorem}[section]

\newtheorem{lem}[thm]{Lemma}

\newtheorem{conj}[thm]{Conjecture}

\numberwithin{equation}{section}

\allowdisplaybreaks

\def\qed{\nopagebreak\hfill{\rule{4pt}{7pt}}}
\def\proof{\noindent {\it{Proof.} \hskip 2pt}}

\begin{document}

\begin{center}

 {\Large \bf Unimodality of partition polynomials related to \\[5pt] Borwein's conjecture}

\end{center}

\begin{center}
 { Janet J. W. Dong}$^{1}$ and {Kathy Q. Ji}$^{2}$  \vskip 2mm

$^{1,2}$ Center for Applied Mathematics,  Tianjin University, Tianjin 300072, P.R. China\\[6pt]
   \vskip 2mm

    $^1$dongjinwei@tju.edu.cn and $^2$kathyji@tju.edu.cn,
\end{center}

\vskip 6mm \noindent {\bf Abstract.}  The objective of this paper is to prove that  the   polynomials $\prod_{k=0}^n(1+q^{3k+1})(1+q^{3k+2})$ are symmetric and unimodal for $n\geq 0$ by an analytical method.



 \vskip 6mm

\section{Introduction}
The study of unimodality of polynomials (or combinatorial sequences) has drawn great attention  in recent decades. There is a remarkable diversity of applicable tools, ranging from analytic to topological, and from representation theory to probabilistic analysis. In this paper, we establish the  unimodality of the   polynomials defined in \eqref{main-eqn} by  refining the method of Odlyzko-Richmond \cite{Odlyzko-Richmond-1982}. Recall that a polynomial
 \[a_0+a_1q+\cdots+a_Nq^n\]
with integer coefficients is called  unimodal if for some $0\leq j\leq N$,
\[a_0\leq a_1\leq \cdots \leq a_j\geq
a_{j+1}\geq \cdots \geq a_N,\]
and is called symmetric if for all $0\leq j\leq N$,
\[a_j=a_{N-j}.\]
See \cite[p. 124, Ex. 50]{Stanley-1997}. It is well-known that the  Gaussian polynomials
\[{n\brack k}=\frac{(1-q^n)(1-q^{n-1})\cdots(1-q^{n-k+1})}{(1-q)(1-q^2)\cdots(1-q^k)}\]
are symmetric and unimodal, as  conjectured by Caylay \cite{Cayley-1856} in 1856  and confirmed by Sylvester \cite{Sylvester-1878} in 1878 based on semi-invariants of binary forms. For more information, we refer to  \cite{Chen-Jia-2022, O'Hara-1990,   Pak-Panova-2013, Proctor-1982}.

 R. C. Entringer  may be the first to investigate the unimodality of  polynomials by an analytical method. By extending the argument of van Lint \cite{Lint-1967}, Entringer \cite{Entringer-1968} showed that the polynomials
 \[(1+q)^2(1+q^2)^2\cdots(1+q^n)^2\]
are unimodal for $n\geq 1$.  This method was greatly extended by Odlyzko and Richmond \cite{Odlyzko-Richmond-1982} to establish the almost unimodality of a class of polynomials of the form
\[(1+q^{a_1})(1+q^{a_2})\cdots(1+q^{a_n})\]
when $n$ is large enough, where $\{a_i\}_{i=1}^\infty$ is a non-decreasing sequence of positive integers.
More precisely, let
\begin{equation}\label{eq1-Odlyzko-Richmond}
    \prod_{i=1}^n(1+q^{a_i})
    =\sum_{m=0}^Nb_n(m)q^m,\qquad \text{where} \quad N=\sum_{i=1}^na_i,
\end{equation}
Odlyzko and Richmond showed that under suitable conditions (conditions (\uppercase\expandafter{\romannumeral1}) and (\uppercase\expandafter{\romannumeral2}) in  Roth and   Szekeres \cite[p. 241]{Roth-Szekeres-1954}) on the infinite sequence $\{a_i\}$, the polynomials \eqref{eq1-Odlyzko-Richmond} are almost unimodal for $n$ sufficiently large, that is, when $n\rightarrow \infty$,
\begin{equation}\label{ineqal-b-g}
 b_n(A) \leq b_n(A+1) \leq \cdots\leq b_n(K)\geq  b_n(K+1)\geq \cdots\geq  b_n(N-A),
\end{equation}
where $A$ is some fixed constant and $K=N/2$ or $K=(N+1)/2$.

When $a_i=i$ for $1\leq i\leq n$ in \eqref{eq1-Odlyzko-Richmond}, Odlyzko and Richmond \cite{Odlyzko-Richmond-1982} verified that the inequality \eqref{ineqal-b-g} holds for $A=1$ when $n\geq 60$. It can be checked that inequality \eqref{ineqal-b-g} also holds for $A=1$ when $n\leq 59$. Hence Odlyzko and Richmond concluded that  the polynomials
\begin{equation}\label{eq2-Odlyzko-Richmond}
 (1+q)(1+q^2)\cdots(1+q^n)
\end{equation}
are unimodal for  $n\geq 1$. The first   proof of the unimodality of the polynomials \eqref{eq2-Odlyzko-Richmond}  was given  by Hughes \cite{Hughes-1977} with the aid of Lie algebra results. Stanley \cite{Stanley-1982}  provided an alternative proof by using the Hard Lefschetz Theorem.  Stanley \cite{Stanley-1980} also established the general result of this type  based on a   result of Dynkin  \cite{Dynkin-1950}.

 When $a_i=2i-1$ for $1\leq i\leq n$ in  \eqref{eq1-Odlyzko-Richmond}, Almkvist \cite{Almkvist-1985} proved that the inequality \eqref{ineqal-b-g} holds for $A=3$ when $n\geq 83$. This leads to  the polynomials
\begin{equation}\label{eq2-Almkvist }
(1+q)(1+q^3)\cdots(1+q^{2n-1})
\end{equation}
are unimodal for $n\geq 27$, except at the coefficient of $q^2$ and $q^{n^2-2}$ conjectured by  Stanley \cite{Stanley-1982}. Pak and Panova \cite{Pak-Panova-2014} showed that the polynomials \eqref{eq2-Almkvist } are strict unimodal  by interpreting
the differences between numbers of certain partitions as Kronecker coefficients of representations of $S_n$.

In \cite{Almkvist-1985}, Almkvist also made the following conjecture.
\begin{conj}[Almkvist] For even $r\geq 2$ or   odd $r\geq 3$ and   $n\geq 11$, the polynomials
\begin{equation}\label{eq-conject}
\prod_{k=1}^n\frac{1-q^{rk}}{1-q^k}
\end{equation}
are unimodal.
\end{conj}

When $r=2$,  the polynomials \eqref{eq-conject} reduces to the polynomials \eqref{eq2-Odlyzko-Richmond}. Almkvist \cite{Almkvist-1987} first showed that the conjecture is true when $r = 4$  by refining the method of Odlyzko-Richmond \cite{Odlyzko-Richmond-1982}. Subsequently,   Almkvist
  \cite{Almkvist-1989} showed that the conjecture is true when  $3\leq r\leq 20$,  $r = 100$ and $101$.

In this paper, we establish the unimodality of the following polynomials.

\begin{thm} \label{main thm}
For $n\geq 0$, the polynomials \begin{equation}\label{main-eqn}
 \prod\limits_{k=0}^{n}(1+q^{3k+1})(1+q^{3k+2})
\end{equation}
are symmetric and unimodal.
\end{thm}

It is worth mentioning that Borwein conjectured that the coefficients of the polynomials
\[ \prod\limits_{k=0}^{n}(1-q^{3k+1})(1-q^{3k+2}) \]
have a repeating sign pattern of $+--$, which has been called as Borwein's conjecture, see Andrews \cite{Andrews-1995}. Recently, Borwein's conjecture has  been  proved by Wang \cite{Wang-2019} by an analytical method.

\section{Preliminaries}

In this section, we collect several identities and inequalities which will be useful in the  proof of Theorem \ref{main thm}.
\begin{align}
    e^{ix}&=\cos(x)+i \sin(x)\label{eq-2.1},\\[5pt]
\cos (2x)&=2\cos^2 (x)-1 \label{eq-2.2}\\[5pt]
&=1-2\sin^2(x)\label{eq-2.2-2},\\[5pt]
\sin (2x)&=2\sin(x) \cos(x) \label{eq-2.2-3},\\[5pt]
2 \sin(\alpha) \cos(\beta) &=  \sin(\alpha+\beta)+\sin(\alpha-\beta)\label{eq-2.3},\\[5pt]
\sin(x)&\geq xe^{- x^2/3}\quad \text{for} \quad 0\leq x\leq 2\label{eq-2.4},\\[5pt]
\cos(x)&\geq e^{-\gamma x^2}\quad \text{for} \quad |x|\leq 1,\,(\gamma =-\log \cos (1)=0.615626....)\label{eq-2.5},\\[5pt]
x-\frac{x^3}{6} &\leq \sin(x)\leq x\quad \text{for} \quad x\geq 0\label{eq-2.6},\\[5pt]
|\cos(x)|&\leq \exp\left(-\frac{1}{2}\sin^2(x)-\frac{1}{4}\sin^4 (x)\right)\quad \text{for} \quad x\geq 0\label{eq-2.7},\\[5pt]
\left|\frac{\sin(nx)}{\sin(x)}\right|&\leq n\label{eq-2.7-1},\\[5pt]
\sum_{k=1}^n\sin^2(kx)&=\frac{n}{2}-\frac{\sin((2n+1)x)}{4\sin(x)}+\frac{1}{4} \label{eq-2.8},\\[5pt]
\sum_{k=1}^n\sin^4(kx)&=\frac{3n}{8}-\frac{\sin((2n+1)x)}{4\sin(x)}+\frac{\sin((2n+1)2x)}{16\sin(2x)}+\frac{3}{16}.\label{eq-2.9}
\end{align}
The identity \eqref{eq-2.1} is Euler's identity, see \cite[p. 4]{Stein-Shakarchi-2003}. For the formulas \eqref{eq-2.2}--\eqref{eq-2.3}  of trigonometric functions, please see \cite[Chapter 8] {Burdette-1973}. The inequalities \eqref{eq-2.4}--\eqref{eq-2.7-1} were proved  by Odlyzko and Richmond \cite[p. 81]{Odlyzko-Richmond-1982}.

It remains to show  \eqref{eq-2.8} and \eqref{eq-2.9}.

\noindent{\it Proofs of \eqref{eq-2.8} and \eqref{eq-2.9}.} First, by \eqref{eq-2.3}, we obtain
\begin{align*}
&2\sin(x)\left(\frac{1}{2}+\sum_{k=1}^n\cos(2kx)\right)\\[5pt]
&\quad =\sin(x)+2\sin(x)\cos(2x)+2\sin(x)\cos(4x)+\cdots+2\sin(x)\cos(2nx)\nonumber\\[5pt]
&\quad \overset{\eqref{eq-2.3}}{=}\sin(x)+\left(\sin(3x)-\sin(x)\right)+\left(\sin(5x)-\sin (3x)\right)\nonumber\\[5pt]
&\quad \quad+\cdots+\left(\sin((2n+1)x)-\sin ((2n-1)x)\right)\nonumber\\[5pt]
&\quad =\sin((2n+1)x).\nonumber
\end{align*}	
Hence, we have
\begin{align}
\sum_{k=1}^n\cos(2kx)
=\frac{\sin((2n+1)x)}{2\sin(x)}-\frac{1}{2}. \label{eq-2.10}
\end{align}	
Using \eqref{eq-2.2-2} and \eqref{eq-2.10}, we deduce that
\begin{align}
\sum_{k=1}^n\sin^2(kx)
&\overset{\eqref{eq-2.2-2}}{=}\frac{n}{2}-\frac{1}{2}\sum_{k=1}^n\cos(2kx)\nonumber\\[5pt]
&\overset{\eqref{eq-2.10}}{=}\frac{n}{2}-\frac{1}{2}\left(\frac{\sin((2n+1)x)}{2\sin(x)}-\frac{1}{2}\right)\nonumber\\[5pt]
&=\frac{n}{2}-\frac{\sin((2n+1)x)}{4\sin(x)}+\frac{1}{4},\nonumber
\end{align}	
which is \eqref{eq-2.8}.

The identity \eqref{eq-2.9} can be derived in the same way. To wit,
\begin{align}
\sum_{k=1}^n\sin^4(kx)
&\overset{\eqref{eq-2.2-2}}{=}\sum_{k=1}^n\left(\frac{1-\cos(2kx)}{2}\right)^2\nonumber\\[5pt]
&\overset{\eqref{eq-2.2}}{=}\frac{3n}{8}-\frac{1}{2}\sum_{k=1}^n\cos(2kx)+\frac{1}{8}\sum_{k=1}^n\cos(4kx)\nonumber\\[5pt]
&\overset{\eqref{eq-2.10}}{=}\frac{3n}{8}-\frac{1}{2}\left(\frac{\sin((2n+1)x)}{2\sin(x)}-\frac{1}{2}\right)+\frac{1}{8}\left(\frac{\sin((2n+1)2x)}{2\sin(2x)}-\frac{1}{2}\right)\nonumber\\[5pt]
&=\frac{3n}{8}-\frac{\sin((2n+1)x)}{4\sin(x)}+\frac{\sin((2n+1)2x)}{16\sin(2x)}+\frac{3}{16},\nonumber
\end{align}	
in agreement with \eqref{eq-2.9}. This completes the proof. \qed

\section{Proof of Theorem \ref{main thm}}

Let $d_n=3(n+1)^2$ and define
\begin{align}\label{defi-gf}
B_n(q)=\prod_{k=0}^{n}(1+q^{3k+1})(1+q^{3k+2})=\sum_{m=0}^{d_n}a_n(m)q^m.
\end{align}

In order to prove Theorem  \ref{main thm}, we first show the following lemma.

\begin{lem} \label{main-lem}
If $n\geq 1$  and $\frac{3n^2}{2}\leq m\leq\frac{3(n+1)^2}{2}$, then
\begin{equation}\label{lem:inequ}
    a_n(m)-a_n(m-1)\geq 0.
\end{equation}
\end{lem}

\proof  We first   show that \eqref{lem:inequ} holds for $n\geq 168$  and $\frac{3n^2}{2}\leq m\leq\frac{3(n+1)^2}{2}$. Putting $q=e^{2i\theta}$ in \eqref{defi-gf},  by \eqref{eq-2.1}, \eqref{eq-2.2} and \eqref{eq-2.2-3},  we derive  that
\begin{align}
	B_n(e^{2i\theta})
	&=\prod_{k=0}^{n}(1+(e^{2i\theta})^{3k+1})(1+(e^{2i\theta})^{3k+2})\nonumber\\[3pt]
	&\overset{\eqref{eq-2.1}}{=}\prod_{k=0}^{n}\left(1+\cos\left(2(3k+1)\theta\right)+i\sin(2(3k+1)\theta)\right)\nonumber\\[3pt]
	&\quad \quad \quad \quad \times \left(1+\cos\left(2(3k+2)\theta\right)+i \sin(2(3k+2)\theta)\right)\nonumber\\[3pt]
	&\overset{\eqref{eq-2.2}\&\eqref{eq-2.2-3}}{=}\prod_{k=0}^{n}\left(2\cos^2((3k+1)\theta)+2i\sin((3k+1)\theta)\cos((3k+1)\theta)\right)\nonumber\\[3pt]
	&\quad\qquad\quad \quad \times\left(2\cos^2((3k+2)\theta)+2i\sin((3k+2)\theta)\cos((3k+2)\theta)\right)\nonumber\\[3pt]
	&\overset{\eqref{eq-2.1}}{=}\prod_{k=0}^{n}4\cos((3k+1)\theta)\cos((3k+2)\theta)\exp(i(3k+1)\theta)\exp(i(3k+2)\theta)\nonumber\\[3pt]
	&=4^{n+1}\exp(id_n\theta)\prod_{k=0}^{n}\cos((3k+1)\theta)\cos((3k+2)\theta).
	\label{defi-ss}
\end{align}
Using Taylor's theorem \cite[p. 47--49]{Stein-Shakarchi-2003}, we find that
\begin{align}
a_n(m)
&=\frac{1}{2\pi i}\int_{-\frac{\pi}{2}}^{\frac{\pi}{2}}\frac{B_n\left(e^{2i\theta}\right)}{\left(e^{2i\theta}\right)^{m+1}}\mathrm{d}\left(e^{2i\theta}\right)\nonumber\\[5pt]
&=\frac{1}{\pi}\int_{-\frac{\pi}{2}}^{\frac{\pi}{2}}B_n\left(e^{2i\theta}\right)e^{-2im\theta}\mathrm{d}\theta\nonumber\\[5pt]
&\overset{\eqref{defi-ss}}{=}\frac{4^{n+1}}{\pi}\int_{-\frac{\pi}{2}}^{\frac{\pi}{2}}\exp({i}(d_n-2m)\theta)\prod_{k=0}^{n}\cos((3k+1)\theta)\cos((3k+2)\theta)\mathrm{d}\theta\nonumber\\[5pt]
&\overset{\eqref{eq-2.1}}{=}\frac{4^{n+1}}{\pi}\int_{-\frac{\pi}{2}}^{\frac{\pi}{2}}\left(\cos((d_n-2m)\theta)+{i}\sin((d_n-2m)\theta)\right)   \prod_{k=0}^{n}\cos((3k+1)\theta)\cos((3k+2)\theta)\mathrm{d}\theta.\nonumber
\end{align}	
Observe that
\[\int_{-\frac{\pi}{2}}^{\frac{\pi}{2}}\sin((d_n-2m)\theta)\prod_{k=0}^{n}\cos((3k+1)\theta)\cos((3k+2)\theta)\mathrm{d}\theta=0,\]
we have therefore,
\[a_n(m)=\frac{2^{2n+3}}{\pi}\int_{0}^{\frac{\pi}{2}}\cos((d_n-2m)\theta)\prod_{k=0}^{n}\cos((3k+1)\theta)\cos((3k+2)\theta)\mathrm{d}\theta.\]

We next show that
\begin{equation}\label{lem:tt2}
\frac{\partial}{\partial m}a_n(m)\geq 0 \quad \text{for} \quad  n\geq 168  \quad \text{and} \quad  \frac{3n^2}{2}\leq m\leq\frac{3(n+1)^2}{2},
\end{equation}
from which, it follows that \eqref{lem:inequ} is valid for $n\geq 168$  and $\frac{3n^2}{2}\leq m\leq\frac{3(n+1)^2}{2}$.

It is easy to see that
\[\frac{\partial}{\partial m}a_n(m)=\frac{2^{2n+4}}{\pi}\int_0^{\frac{\pi}{2}}\theta\sin\left((d_n-2m)\theta\right)\prod_{k=0}^{n}\cos((3k+1)\theta)\cos((3k+2)\theta)\mathrm{d}\theta.\]
Let $d_n-2m=\mu$, and let
\[I_n(\mu)=\int_0^{\frac{\pi}{2}}\theta\sin\left(\mu\theta\right)\prod_{k=0}^{n}\cos((3k+1)\theta)\cos((3k+2)\theta)\mathrm{d}\theta.\]
Under the condition that $\frac{3n^2}{2}\leq m\leq\frac{3(n+1)^2}{2}$, we see that
\begin{equation}\label{defi-mu}
   0\leq \mu=d_n-2m\leq 6n+3.
\end{equation}
To prove \eqref{lem:tt2}, it suffices  to show that
\begin{equation}\label{lem:tt2-i}
I_n(\mu)\geq 0 \quad \text{for} \quad  n\geq 168  \quad \text{and} \quad  0\leq \mu \leq 6n+3.
\end{equation}
To this end,  we  write
\begin{align}
	I_n(\mu)
	&=\left \{\int_0^{\frac{\pi}{6n+4}}+\int_{\frac{\pi}{6n+4}}^{\frac{\pi}{2}}\right \}\theta\sin\left(\mu\theta\right)\prod_{k=0}^{n}\cos((3k+1)\theta)\cos((3k+2)\theta)\mathrm{d}\theta\nonumber\\[5pt]
	&=I^{(1)}_n(\mu)+I^{(2)}_n(\mu).\nonumber
\end{align}	
We next show that
 \begin{equation}\label{lem:mc}
 I^{(1)}_n(\mu)\geq |I^{(2)}_n(\mu)| \quad \text{for} \quad  n\geq 168  \quad \text{and} \quad  0\leq \mu \leq 6n+3,
 \end{equation}
 which implies \eqref{lem:tt2-i}.

We first evaluate the value of  $ I^{(1)}_n(\mu)$, which is defined by
\begin{equation}\label{defi-i1}
  I^{(1)}_n(\mu):=\int_0^{\frac{\pi}{6n+4}} \theta\sin\left(\mu\theta\right)\prod_{k=0}^{n}\cos((3k+1)\theta)\cos((3k+2)\theta)\mathrm{d}\theta.
\end{equation}
When $0\leq \theta\leq \frac{1}{3n+2}$, by \eqref{defi-mu},  we have
\[0\leq \mu \theta \leq 2\ \text{ and }\ 0\leq(3k+1)\theta\leq (3k+2)\theta \leq 1\  \text{ for }\  0\leq k\leq n,\]
so that
\begin{align}
&\theta\sin\left(\mu\theta\right)\prod_{k=0}^{n}\cos((3k+1)\theta)\cos((3k+2)\theta)\nonumber\\[3pt]
&\overset{\eqref{eq-2.4} \& \eqref{eq-2.5}}{\geq} \mu\theta^2\exp\left(-\frac{\mu^2\theta^2}{3}\right)\exp\left(-\gamma\theta^2\sum_{k=0}^n\left((3k+1)^2+(3k+2)^2\right)\right)\nonumber\\[3pt]
&\geq \mu\theta^2\exp\left(-\frac{(6n+3)^2\theta^2}{3}\right)\exp\left(-\gamma\theta^2\left(6n^3+18n^2+17n+5\right)\right)\nonumber \quad (\text{by }  0\leq \mu\leq 6n+3)\\[3pt]
&=\mu\theta^2\exp\left(-\theta^2n^3\left(\left(\frac{12}{n}+\frac{12}{n^2}+\frac{3}{n^3}\right)+\gamma\left(6+\frac{18}{n}+\frac{17}{n^2}+\frac{5}{n^3}\right)\right)\right)\nonumber \\[3pt]
&\geq\mu\theta^2\exp\left(-cn^3\theta^2\right) \quad (\text{by }  n\geq 168),\label{eva-inter}
\end{align}	
where $c=3.832$.
Applying \eqref{eva-inter} to \eqref{defi-i1}, we find that  when $n\geq 168$ and $0\leq \mu \leq 6n+3$,
\begin{align}
I^{(1)}_n(\mu)&=\int_0^{\frac{\pi}{6n+4}}\theta\sin\left(\mu\theta\right)\prod_{k=0}^{n}\cos((3k+1)\theta)\cos((3k+2)\theta)\mathrm{d}\theta\nonumber\\[3pt]
&\geq \int_0^{\frac{1}{3n+2}}\theta\sin\left(\mu\theta\right)\prod_{k=0}^{n}\cos((3k+1)\theta)\cos((3k+2)\theta)\mathrm{d}\theta\nonumber\\[3pt]
&\geq \int_0^{\frac{1}{3n+2}}\mu\theta^2\exp\left(-cn^3\theta^2\right)\mathrm{d}\theta\nonumber\\[3pt]
&=\left \{\int_0^{\infty}-\int_{\frac{1}{3n+2}}^{\infty}\right \}\mu\theta^2\exp\left(-cn^3\theta^2\right)\mathrm{d}\theta\nonumber\\[3pt]
&=\frac{\mu}{2c^{\frac{3}{2}}n^{\frac{9}{2}}}\left(\int_0^{\infty}v^{\frac{1}{2}}e^{-v}\mathrm{d}v-\int_{\frac{cn^3}{(3n+2)^2}}^{\infty}v^{\frac{1}{2}}e^{-v}\mathrm{d}v\right)\nonumber\\[3pt]
&=\frac{\mu}{2c^{\frac{3}{2}}n^{\frac{9}{2}}}\left(\frac{\sqrt{\pi}}{2}-\int_{\frac{cn^3}{(3n+2)^2}}^{\infty}v^{\frac{1}{2}}e      ^{-v}\mathrm{d}v\right).\nonumber
\end{align}	
Observe that  when $n\geq 168$,
\[\frac{c n^3}{(3n+2)^2}\geq \frac{c\cdot168^3}{(3\times 168+2)^2},\]
so
\[\int_{\frac{cn^3}{(3n+2)^2}}^{\infty}v^{\frac{1}{2}}e^{-v}\mathrm{d}v\leq\int_{\frac{c\cdot168^3}{(3\times168+2)^2}}^{\infty}v^{\frac{1}{2}}e^{-v}\mathrm{d}v\leq 1.29\times 10^{-30}.\]
Consequently, when  $n\geq 168$ and $0\leq \mu \leq 6n+3$,
\begin{align}
I^{(1)}_n(\mu)&\geq\frac{\frac{\sqrt{\pi}}{2}-1.29\times10^{-30}}{2\times3.832^{\frac{3}{2}}}\cdot \frac{\mu}{n^{\frac{9}{2}}}\geq\frac{0.8862\mu}{15.2 n^{\frac{9}{2}}}\geq \frac{0.0583\mu}{n^{\frac{9}{2}}}.\label{eq-I1}
\end{align}	

We now turn to estimate the value of  $ I^{(2)}_n(\mu)$, which is defined by
\begin{equation}\label{defi-i2}
  I^{(2)}_n(\mu)=
	 \int_{\frac{\pi}{6n+4}}^{\frac{\pi}{2}}\theta\sin\left(\mu\theta\right)\prod_{k=0}^{n}\cos((3k+1)\theta)\cos((3k+2)\theta)\mathrm{d}\theta.
\end{equation}
When $\frac{\pi}{6n+4}\leq \theta \leq \frac{\pi}{2}$, by \eqref{eq-2.7}, \eqref{eq-2.8} and \eqref{eq-2.9}, we deduce that
\begin{align}
&\left|\prod_{k=0}^{n}\cos((3k+1)\theta)\cos((3k+2)\theta)\right|\nonumber\\[5pt]
&\overset{\eqref{eq-2.7}}{\leq}\exp\left (-\frac{1}{2}\sum_{k=0}^{n}\left(\sin^2((3k+1)\theta)+\sin^2((3k+2)\theta)\right)\right. \nonumber\\[5pt]
&\quad \quad\quad\quad \quad\quad\quad\quad\quad\quad\left. -\frac{1}{4}\sum_{k=0}^{n}\left(\sin^4((3k+1)\theta)+\sin^4((3k+2)\theta)\right)\right)	\nonumber\\[8pt]
&=\exp\left (-\frac{1}{2}\left(\sum_{k=1}^{3n+2}\sin^2(k\theta)-\sum_{k=1}^{n}\sin^2(3k\theta)\right)-\frac{1}{4}\left(\sum_{k=1}^{3n+2}\sin^4(k\theta)-\sum_{k=1}^{n}\sin^4(3k\theta)\right)\right )\nonumber\\[10pt]
&\overset{\eqref{eq-2.8}\&\eqref{eq-2.9}}{=}\exp\left (-\frac{11(n+1)}{16}+\frac{3\sin((6n+5)\theta)}{16\sin(\theta)}-\frac{\sin((6n+5)2\theta)}{64\sin(2\theta)}\right.\nonumber\\[5pt]
&\quad \quad\quad\quad \quad\quad\quad \quad\quad\quad\left.  -\frac{3\sin((2n+1)3\theta)}{16\sin(3\theta)}+\frac{\sin((2n+1)6\theta)}{64\sin(6\theta)}\right):=E(n).\nonumber
\end{align}	
We proceed to prove that
\begin{equation}\label{eval-value-e}
	E(n)<\exp\left(-0.163n-0.031\right) \quad \text{for} \quad  \frac{\pi}{6n+4}\leq \theta \leq \frac{\pi}{2}  \quad \text{and} \quad  n\geq 168.
\end{equation}
The proof of \eqref{eval-value-e} is divided into two steps. When $\frac{\pi}{6n+4}\leq \theta \leq \frac{\pi}{6}$, using \eqref{eq-2.6} and \eqref{eq-2.7-1},  we obtain
\begin{align}
	E(n)
	&\leq\exp\left (-\frac{11(n+1)}{16}+\frac{3}{16\sin(\theta)}+\frac{1}{64\sin(2\theta)}+\frac{3}{16\sin(3\theta)}+\left|\frac{\sin((2n+1)6\theta)}{64\sin(6\theta)}\right|\right )\nonumber\\[5pt]	
	&\overset{\eqref{eq-2.6}\&\eqref{eq-2.7-1}}{\leq}\exp\left (-\frac{11(n+1)}{16}+\frac{3}{16\left(\frac{\pi}{6n+4}\left(1-\frac{\left(\frac{\pi}{6n+4}\right)^2}{6}\right)\right)}+\frac{1}{64\left(\frac{\pi}{3n+2}\left(1-\frac{\left(\frac{\pi}{3n+2}\right)^2}{6}\right)\right)}\right.\nonumber\\[5pt]
	& \left.\qquad\quad+\frac{3}{16\left(\frac{3\pi}{6n+4}\left(1-\frac{\left(\frac{3\pi}{6n+4}\right)^2}{6}\right)\right)}+\frac{2n+1}{64}\right ) \quad \left(\text{by } \frac{\pi}{6n+4}\leq \theta \leq \frac{\pi}{6}\right).\label{eval-e}
\end{align}	
Applying
\[1-\frac{\left(\frac{\pi}{6n+4}\right)^2}{6}\geq1-\frac{\left(\frac{\pi}{3n+2}\right)^2}{6}\geq1-\frac{\left(\frac{3\pi}{6n+4}\right)^2}{6}\]
to \eqref{eval-e}, we derive that
\begin{align}
E(n)&\leq\exp\left (-\frac{42n+43}{64}+\frac{3}{16\left(\frac{\pi}{6n+4}\left(1-\frac{\left(\frac{3\pi}{6n+4}\right)^2}{6}\right)\right)}+\frac{1}{64\left(\frac{\pi}{3n+2}\left(1-\frac{\left(\frac{3\pi}{6n+4}\right)^2}{6}\right)\right)}\right.\nonumber\\[5pt]
&\phantom{=\;\;}\left.~~~~~~~~~~~~~\qquad\quad+\frac{3}{16\left(\frac{3\pi}{6n+4}\left(1-\frac{\left(\frac{3\pi}{6n+4}\right)^2}{6}\right)\right)}\right )\nonumber\\[5pt]
&=\exp\left (-\frac{42n+43}{64}+\frac{33}{128\left(\frac{\pi}{6n+4}\left(1-\frac{\left(\frac{3\pi}{6n+4}\right)^2}{6}\right)\right)}\right )\nonumber\\[5pt]
&=\exp\left (-\frac{42n+43}{64}+\frac{33(6n+4)}{128\pi\left(1-\frac{6\pi^2}{(12n+8)^2}\right)}\right ).\nonumber
\end{align}
Note that when $n\geq 168$,
\[1-\frac{6\pi^2}{(12n+8)^2}\geq1-\frac{6\pi^2}{(12\times168+8)^2}=1-\frac{3\pi^2}{2048288},\]
so when $ \frac{\pi}{6n+4}\leq \theta \leq \frac{\pi}{6}$ and $n\geq 168$,
\begin{align}
E(n)&\leq\exp\left (-\frac{42n+43}{64}+\frac{33(6n+4)}{128\pi\left(1-\frac{3\pi^2}{2048288}\right)}\right )\nonumber\\[5pt]
&=\exp\left (\left(-\frac{21}{32}+\frac{99}{64\pi\left(1-\frac{3\pi^2}{2048288}\right)}\right)n-\frac{43}{64}+\frac{33}{32\pi\left(1-\frac{3\pi^2}{2048288}\right)}\right )\nonumber\\[5pt]
&<\exp\left (-0.163n-0.343\right ).
\label{eq-E-1}
\end{align}
When $\frac{\pi}{6}\leq \theta \leq \frac{\pi}{2}$, by \eqref{eq-2.7-1}, we deduce that
\begin{align}
	E(n)
	&\leq\exp\left (-\frac{11(n+1)}{16}+\frac{3}{16\sin(\theta)}+\left|\frac{\sin((6n+5)2\theta)}{64\sin(2\theta)}\right|+\left|\frac{3\sin((2n+1)3\theta)}{16\sin(3\theta)}\right|\right.\nonumber\\[5pt]
	&\quad \quad\quad \quad\quad \quad\left. +\left|\frac{\sin((2n+1)6\theta)}{64\sin(6\theta)}\right|\right )\nonumber\\[8pt]	
	&\overset{\eqref{eq-2.7-1}}{\leq}\exp\left (-\frac{11(n+1)}{16}+\frac{3}{16\sin(\frac{\pi}{6})}+\frac{6n+5}{64}+\frac{3(2n+1)}{16}+\frac{2n+1}{64}\right )\nonumber\\[5pt]
	&=\exp\left(-\frac{3}{16}n-\frac{1}{32}\right)<\exp\left(-0.187n-0.031\right).
	\label{eq-E-2}
\end{align}	
Combining \eqref{eq-E-1} and \eqref{eq-E-2} yields \eqref{eval-value-e}. Applying \eqref{eval-value-e} to  \eqref{defi-i2}, and in view of \eqref{eq-2.6} and \eqref{eq-I1}, we derive that when $n\geq 168$,
\begin{align}
|I^{(2)}_n(\mu)|&\overset{\eqref{eq-2.6}}{<} \mu \exp\left (-0.163n-0.031\right )\int_{\frac{\pi}{6n+4}}^{\frac{\pi}{2}}\theta^2\mathrm{d}\theta	\nonumber\\[5pt]
&\leq \frac{\mu\pi^3}{3}\left(\frac{1}{8}-\frac{1}{(6n+4)^3}\right) \exp\left (-0.163n-0.031\right )	\nonumber\\[5pt]
&= \frac{\mu\pi^3}{3}\left(\frac{1}{2}-\frac{1}{6n+4}\right)\left(\frac{1}{2^2}+\frac{1}{2(6n+4)}+\frac{1}{(6n+4)^2}\right) \exp\left (-0.163n-0.031\right )	\nonumber\\[5pt]
&< \frac{\mu\pi^3}{3}\cdot\frac{3}{4}\cdot\left(\frac{1}{2}-\frac{1}{6n+4}\right)\exp\left (-0.163n-0.031\right )	\nonumber\\[5pt]
&\overset{\eqref{eq-I1}}{\leq}\frac{\pi^3n^{\frac{9}{2}}}{4\times 0.0583}\left(\frac{1}{2}-\frac{1}{6n+4}\right)\exp\left (-0.163n-0.031\right )I^{(1)}_n(\mu).\nonumber
\end{align}	
Define
\[f(n):=\frac{\pi^3n^{\frac{9}{2}}}{4\times 0.0583}\left(\frac{1}{2}-\frac{1}{6n+4}\right)\exp\left (-0.163n-0.031\right ).\]
To show \eqref{lem:mc}, it remains to show that $f(n)<1$ for $n\geq 168$. We claim that $f'(n)<0$ for $n\geq 168$. Since  $f(n)>0$ for $n\geq 168$, we have
\begin{equation}\label{lem:temp2}
 \frac{\mathrm{d}}{\mathrm{d}n} {f(n)}=\frac{\mathrm{d}}{\mathrm{d}n} e^{\ln f(n)}=f(n)\frac{\mathrm{d}}{\mathrm{d}n} \ln{f(n)}.
\end{equation}
Observe that when $n\geq 168$,
\begin{align}
\frac{\mathrm{d}}{\mathrm{d}n} \ln{f(n)}
&=\frac{9}{2n}+\frac{6}{(3n+1)(6n+4)}-0.163\nonumber\\[5pt]
&\leq\frac{9}{2\times168}+\frac{6}{(3\times168+1)(6\times168+4)}-0.163< -0.13<0.\nonumber
\end{align}
Hence, we derive from  \eqref{lem:temp2} that $f'(n)<0$ for $n\geq 168$, and the claim is proved.  Consequently,  $f(n)\leq f(168)<0.851$ when $n\geq 168$.
Therefore, \eqref{lem:mc} is valid, and so \eqref{lem:tt2} is valid. This leads to \eqref{lem:inequ} holds for $n\geq 168$  and $\frac{3n^2}{2}\leq m\leq\frac{3(n+1)^2}{2}$. Using Maple,
we can check that \eqref{lem:inequ}  also holds  for $n< 168$  and $\frac{3n^2}{2}\leq m\leq\frac{3(n+1)^2}{2}$. Thus the lemma is proved. \qed

We conclude this paper with the proof of Theorem \ref{main thm}.

\noindent{\it Proof of Theorem \ref{main thm}.} When $n\geq 0$, we first show that $B_n(q)$ is a symmetric polynomial. Replacing $q$ by $q^{-1}$ in \eqref{defi-gf}, we deduce that
\begin{align*}
B_n(q^{-1})&=\prod_{k=0}^{n}(1+q^{-(3k+1)})(1+q^{-(3k+2)}) \\[5pt]
&=q^{-d_n}\prod_{k=0}^{n}(1+q^{(3k+1)})(1+q^{(3k+2)}) \\[5pt]
&=q^{-d_n} B_n(q).
\end{align*}
To wit,
\[B_n(q)=q^{d_n}B_n(q^{-1}),\]
from which,  it follows that  $B_n(q)$ is   symmetric.

We proceed to show that the polynomial  $B_n(q)$ is unimodal  by induction on $n$. When $n=0$, we have
\[B_0(q)=(1+q)(1+q^2)=1+q+q^2+q^3.\]
Clearly, the coefficients of $B_0(q)$ are unimodal.

Suppose that  $B_{n-1}(q)$ is  unimodal for $n\geq 1$, namely, for  $n\geq 1$ and $1\leq m\leq \lfloor \frac{d_{n-1}}{2}\rfloor$,
\begin{equation}\label{ineq-n-1}
a_{n-1}(m)\geq a_{n-1}(m-1).
\end{equation}
We intend to show that $B_n(q)$ is unimodal. Since $B_n(q)$ is a symmetric polynomial, it suffices to show that  for $n\geq 1$ and $1\leq m\leq \lfloor \frac{d_{n}}{2}\rfloor$,
\begin{equation}\label{ineq-n}
a_{n}(m)\geq a_{n}(m-1).
\end{equation}
Observe that
\[B_n(q)=\left(1+q^{3n+1}\right)\left(1+q^{3n+2}\right)B_{n-1}(q),\]
which implies the following recurrence relation:
\begin{align}\label{thm:recurrence}
a_n(m)=a_{n-1}(m)+a_{n-1}(m-3n-1)+a_{n-1}(m-3n-2)+a_{n-1}(m-6n-3).
\end{align}	
It's evident from \eqref{ineq-n-1} and \eqref{thm:recurrence}  that \eqref{ineq-n} holds for $n\geq 1$ and $1\leq m\leq \lfloor \frac{d_{n-1}}{2}\rfloor$.
In view of Lemma \ref{main-lem}, we see that \eqref{ineq-n} also holds for $n\geq 1$ and $\lceil \frac{d_{n-1}}{2}\rceil\leq m\leq \lfloor \frac{d_{n}}{2}\rfloor$.
Hence, we conclude that \eqref{ineq-n} is valid for $n\geq 1$ and $1\leq m\leq \lfloor \frac{d_{n}}{2}\rfloor$, and so  $B_n(q)$ is unimodal. Thus, we complete the proof of Theorem \ref{main thm}. \qed
 \vskip 0.2cm
\noindent{\bf Acknowledgment.} This work
was supported by   the National Science Foundation of China.

\end{document}